\documentclass[11pt, a4paper]{amsart}
\usepackage{fullpage}
\usepackage{amscd}
\usepackage{mathtools}
\usepackage{amsmath}
\usepackage{amsthm}
\usepackage{amssymb}
\usepackage{graphicx}
\usepackage[colorlinks,linkcolor=blue,filecolor=blue,citecolor=blue,urlcolor=blue,bookmarks=false,hypertexnames]{hyperref}
\usepackage{tikz}
\usepackage{caption}
\usepackage{caption,subcaption}
\usepackage{paralist}
\newcommand\naturals{\mathbb{N}}

\def\opn#1#2{\def#1{\operatorname{#2}}} % to make operators
\let\iso=\simeq
\def\opn#1#2{\def#1{\operatorname{#2}}} % to make operators
\opn\width{width}
\opn\rank{rank} 
\opn\min{min}
\opn\max{max}
\opn\height{height}
\opn\link{link}
\opn\id{id}
\makeatletter
\@addtoreset{equation}{section}
\def\theequation{\thesection.\@arabic \c@equation}
\makeatother

\theoremstyle{plain}

\newtheorem{theorem}{Theorem}

\newtheorem{Lemma}[equation]{Lemma}
\newtheorem{proposition}[equation]{Proposition}

\newtheorem{setup}[equation]{Setup}

\theoremstyle{definition}

\newtheorem{definition}[equation]{Definition}
\newenvironment{definitionbox}[1][]{%
\begin{definition}[#1]\pushQED{\qed}}{\popQED \end{definition}}

\newtheorem{example}[equation]{Example}
\newenvironment{examplebox}[1][]{%
\begin{example}[#1]\pushQED{\qed}}{\popQED \end{example}}

\newtheorem{discussion}[equation]{Discussion}

\newtheorem{Notation}[equation]{Notation}

\newtheorem{observation}[equation]{Observation}
\newenvironment{observationbox}[1][]{%
\begin{observation}[#1]\pushQED{\qed}}{\popQED \end{observation}}

\title {On Cohen-Macaulay posets of dimension two and permutation graphs}

	\author{Rizwan Jahangir}
\address{Sabanci University, Faculty of Engineering and Natural Sciences, Orta Mahalle, Tuzla 34956, Istanbul, Turkey}
\email{rizwan@sabanciuniv.edu}

\author{Dharm Veer}
\address{Indian Institute of Technology Gandhinagar, Palaj, Gujarat 382355. India
    \newline
Chennai Mathematical Institute, Siruseri, Tamilnadu 603103. India}
\email{dharm.v@iitgn.ac.in}

\thanks{DV was supported by a grant of IIT Gandhinagar and was partly
supported by an Infosys Foundation fellowship.}



\begin{document}

\begin{abstract}
We characterize Cohen-Macaulay posets of dimension two;
they are precisely the shellable and strongly connected posets of dimension two.
We also give a combinatorial description of these posets.
Using the fact that co-comparability graph of a 2-dimensional poset is a permutation graph, 
we characterize Cohen-Macaulay permutation graphs.
\end{abstract}
\maketitle
\section{Introduction}

The notion of Cohen-Macaulay posets was first defined in Baclawski's thesis~\cite{baclawskithesis} purely combinatorially (see \cite[Section\ 3]{baclawskicmdefinition80}). 
Reisner~\cite{ReisnerCMrings} and Stanley~\cite{stanleycmrings75} independently provided a ring theoretic definition of Cohen-Macaulay posets, with Reisner~\cite[Theorem\ 1]{ReisnerCMrings} demonstrating the equivalence of the two definitions.
Stanley~\cite{stanleyupperboundconj75} used Reisner's result to prove the upper bound conjecture for spheres.
This was the inception of the Stanley-Reisner theory, which served as a bridge connecting the fields of combinatorics, commutative algebra, and topology.

The purpose of this article is to characterize the Cohen-Macaulay posets of dimension two, indeed we show that they are shellable and strongly connected. 
A {\em linear extension} $\pi$ of a poset $P$ is a linear order on the underlying set of $P$ such that $x \leq y$ in $\pi$ whenever $x \leq y$ in $P$.
A poset $P$ is an intersection of a family of linear extensions $\pi_1,\ldots,\pi_d$ if $x\leq y$ in $P$ if and only if  $x\leq y$ in $\pi_i$, for all $1\leq i\leq d.$ 
The {\em dimension} of a poset $P$ is the least integer $d$ such that $P$ can be expressed as the intersection of $d$ linear extensions of $P$. 
The dimension of a poset was defined by Dushnik and Miller~\cite{millerposetdimension}.

A poset is said to be Cohen-Macaulay over a field $K$ if its order complex is Cohen-Macaulay over $K$ (see Section~\ref{sec:pre} for definitions).
Similarly, we say that a poset is shellable (resp.~strongly connected) if its order complex is shellable (resp.~strongly connected).
It is known that shellable posets are Cohen-Macaulay over any field~\cite[Theorem\ 5.1.13]{BH93}, and 
Cohen-Macaulay posets are strongly connected (see \cite[Proposition\ 11.7]{Bjornertopologicalmethods}).

Let $P$ be a poset. 
We say that $P$ is an {\em antichain} if any two distinct elements of $P$ are incomparable.
For $p\in P$, {\em height} of $p$ is the rank of the induced subposet of $P$ which consists of all $q\in P$ with $q\leq p$.
In this paper, we prove the following:

\begin{theorem}\label{thm:cmpermutation}
    Let $P$ be a finite poset of dimension two.
    Then the following are equivalent:
   \begin{enumerate}
       \item\label{thm:cmpermutation:shellable}
        $P$ is shellable.
       
       \item \label{thm:cmpermutation:cm}
        $P$ is Cohen-Macaulay.

       \item \label{thm:cmpermutation:stronglycon}
           $P$ is strongly connected.

       \item \label{thm:cmpermutation:poset}
        $P$ is an antichain or $P$ is pure and the induced subposet of $P$ consisting of height $i$ and height $i+1$ elements is connected for all $0 \leq i \leq rank(P) -1$.  
   \end{enumerate}
\end{theorem}

We use Theorem~\ref{thm:cmpermutation} to characterize Cohen-Macaulay permutation graphs.
The proof uses the fact~\cite[Theorem\ 1]{GRU83comparabilitygraphs} that a permutation graph is a co-comparability graph of a poset of dimension at most two,
and that the Stanley-Reisner ideal of the order complex of a poset coincides with the edge ideal of its co-comparability graph.
Section~\ref{sec:pre} contains the definitions and preliminaries. The proof of the theorem is given in Section~\ref{sec:proof:thm:cmpermutation}.

\subsection*{Acknowledgements}
This project was started when the second author visited Ayesha Asloob Qureshi at Sabancı University, Turkey, he thank her for the hospitality. 
Both authors thank her for several helpful discussions.
\section{Preliminaries}\label{sec:pre}
Throughout this article, all posets are finite  and all graphs are simple and finite.

Let $K$ be field and let $\Delta$ be a finite simplicial complex. 
For a face $\sigma$ of $\Delta$, define the {\em link} of $\sigma$ in $\Delta$, denoted by $\link(\Delta,\sigma)$ to be the subcomplex 
$\{\tau\in \Delta : \tau\cup \sigma\in \Delta, \tau\cap\sigma=\varnothing\}$.
The simplicial complex $\Delta$ is said to be {\em Cohen-Macaulay} over $K$ if $\widetilde{H}_i(\link(\Delta,\sigma), K)= 0$ for all $i< \dim(\link(\Delta,\sigma))$ for every face $\sigma$ of $\Delta$.
Here, $\widetilde{H}_i(\_, K)$ is the $i$-th reduced homology group with coefficients in $K$.
Reisner~\cite[Theorem\ 1]{ReisnerCMrings} proved that $\Delta$ is Cohen-Macaulay over $K$ if and only if the Stanley-Reisner ring associated to $\Delta$ (over $K$) is Cohen-Macaulay.

Let $P$ be a poset.
For $x,y \in P$, we say that {\em y covers x}, denoted by $x\lessdot y$, if $x<y$ and there is no $z\in P$ with $x<z<y$.
A {\em chain} $C$ of $P$ is a totally ordered subset of $P$. 
The {\em length} of a chain $C$ of $P$ is $\#C - 1$, where $\#C$ is the cardinality of $C$. 
The {\em rank} of $P$, denoted by $\rank(P)$, is the maximum of the lengths of chains in $P$.
A poset is called {\em pure} if all maximal chains of $P$ have the same length.
An {\em induced subposet} $Q$ of $P$ is a poset on a subset of the underlying set $P$ such that for every $x,y\in Q$, $x\leq y$ in $Q$ if and only if $x\leq y$ in $P$.
The order complex $\Delta(P)$ of $P$ is a simplicial complex on the underlying set of $P$ whose faces are chains of $P$.

 Let $P$ and $Q$ be two posets on disjoint sets. 
 The {\em disjoint union} of posets $P$ and $Q$ is the poset $P+Q$ on the set $P\cup Q$ with the following order: 
 if $x,y \in P+Q$, then $x\leq y$ if either $x, y \in P$ and $x \leq y$ in $P$ or $x, y \in Q$ and $x \leq y$ in $Q$. 
 A poset which can be written as disjoint union of two posets is called {\em disconnected}; otherwise the poset is called {\em connected}.

Let $P$ be a pure poset and $\Delta(P)$ be its order complex. We say that $\Delta(P)$ is {\em strongly connected} if for any two maximal chains $\gamma$ and $\gamma'$ of $P$, 
there is a sequence $\sigma_0,\sigma_1, \ldots, \sigma_k$ of maximal chains of $P$ such that $\sigma_0=\gamma$, $\sigma_k=\gamma'$, and $\sigma_i\cap\sigma_{i+1}$ is a chain of length $\rank(P) -1$.
It is known that Cohen-Macaulay complexes are strongly connected (see \cite[Proposition\ 11.7]{Bjornertopologicalmethods}).

The order complex $\Delta(P)$ is called {\em shellable} if the maximal chains of $P$ admit a linear order $\gamma_0,\ldots, \gamma_m$ such that 
for all $1\leq j<i \leq m$, there exists a $v\in \gamma_i \setminus \gamma_j$ and some $k\in [i-1]$ with $\gamma_i\setminus \gamma_k =\{v\}$.
A linear order satisfying the definition is called a {\em shelling order} on $P$.
A shellable complex is Cohen-Macaulay over any field~\cite[Theorem\ 5.1.13]{BH93}.

The {\em co-comparability graph} $G$ of a poset $P$ is a graph on the underlying set of $P$ such that $\{x,y\}$ is an edge of $G$ if and only if $x$ and $y$ are incomparable in $P$.
We say that a graph is a co-comparability graph if it is a co-comparability graph of some poset.

\begin{figure}[h]
\begin{center}
    \begin{tikzpicture}[scale=1.3]
        \draw[thick] (.5,0)--(4.5,0)
                     (.5,1)--(4.5,1)
                     (.5,2)--(4.5,2);
          
        \draw[]      (1,2)--(3,1)--(4,0)
        (2,2)--(1,1)--(2,0)
        (3,2)--(2,1)--(1,0)
        (4,2)--(4,1)--(3,0) 
;

\filldraw[black] (4.6,2) circle (0pt) node[anchor=west]  {$l_0$};
\filldraw[black] (4.6,1) circle (0pt) node[anchor=west]  {$l_1$};
\filldraw[black] (4.6,0) circle (0pt) node[anchor=west]  {$l_2$};

\filldraw[black] (1,0) circle (.5pt) node[anchor=north] {2};
\filldraw[black] (2,0) circle (.5pt) node[anchor=north] {1};
\filldraw[black] (3,0) circle (.5pt) node[anchor=north] {4};
\filldraw[black] (4,0) circle (.5pt) node[anchor=north]  {3};

\filldraw[black] (1,1) circle (.5pt) node[anchor=north] {1};
\filldraw[black] (2,1) circle (.5pt) node[anchor=north] {2};
\filldraw[black] (3,1) circle (.5pt) node[anchor=north] {3};
\filldraw[black] (4,1) circle (.5pt) node[anchor=north]  {4};

\filldraw[black] (1,2) circle (.5pt) node[anchor=south] {3};
\filldraw[black] (2,2) circle (.5pt) node[anchor=south] {1};
\filldraw[black] (3,2) circle (.5pt) node[anchor=south] {2};
\filldraw[black] (4,2) circle (.5pt) node[anchor=south]  {4};

\end{tikzpicture}
\caption{A concatenation of $2$ permutation diagrams}\label{figure:permutation}
\end{center}
\end{figure}
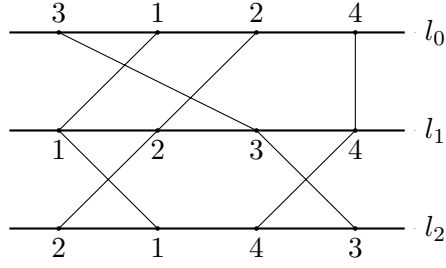

Let $l_0,\ldots, l_k$ be horizontal lines each labeled from left to right by permutations of $[n] = \{1,\ldots,n\}$.
For each $i\in [n]$, the curve $f_i$ consists of $k$ straight line segments which join $i$ on $l_r$ to $i$ in $l_{r+1}$, for $0\leq r \leq k-1$.
When $k=1$, such a diagram is called a {\em permutation diagram}.
When $k\geq 2$, it is called concatenation of $k$ permutation diagrams.
Figure~\ref{figure:permutation} gives an example of a concatenation of $2$ permutation diagrams for $n=4$.

The {\em intersection graph} $G$ of the concatenation of $k$ permutation diagrams is a graph on $[n]$
such that $\{i,j\}$ is an edge of $G$ if and only if $f_i$ intersects with $f_j$. 
It was shown by Golumbic et al.~\cite[Theorem\ 1]{GRU83comparabilitygraphs} that a graph is a co-comparability graphs if and only if it is an
intersection graph of concatenation of $k$ permutation diagrams.

A graph $G$ is called a {\em permutation graph} if it is the intersection graph of a permutation diagram (i.e., $k=1$). 
Observe that the dimension of a poset is at most two if and only if its co-comparability graph is a permutation graph.
Also, it follows from the definition that the dimension of a poset is one if and only if it is a linear order.

\section{Proof of the Theorem~\ref{thm:cmpermutation}}\label{sec:proof:thm:cmpermutation}

For \eqref{thm:cmpermutation:shellable} $\implies$ \eqref{thm:cmpermutation:cm} see~\cite[Theorem\ 5.1.13]{BH93}, and \eqref{thm:cmpermutation:cm}  $\implies$ \eqref{thm:cmpermutation:stronglycon} follows from~\cite[Proposition\ 11.7]{Bjornertopologicalmethods}.
For~\eqref{thm:cmpermutation:stronglycon} $\implies$ \eqref{thm:cmpermutation:poset} of the theorem, we prove a more general statement in the following lemma.

\begin{Lemma} \label{lem:stronglyconnected}
   Let $P$ be a strongly connected poset. 
   Then $P$ is an antichain or the induced subposet of $P$ consisting of height $i$ and height $i+1$ elements is connected for all $0 \leq i \leq \rank(P) -1$.
\end{Lemma}

\begin{proof}
    Clearly, antichains are strongly connected. So we may assume that $\rank(P)\geq 1$.
    We proceed by contradiction. 
    Fix an $i$ with $0 \leq i \leq \rank(P) -1$ such that the induced subposet $Q$ of $P$ consisting of height $i$ and height $i+1$ elements is disconnected. 
    Assume that $Q$ is the disjoint union of two subposets $Q_1$ and $Q_2$.
    Since $P$ is pure, every maximal chain of $P$ contains an element of height $i$ and an element of height $i+1$. 
    Thus, $Q_j$ contains at least one element of height $i$ and at least one element of height $i+1$ for all $j=1,2$.
    
    Let $m_1$ and $m_2$ be two maximal chains of $P$ such that $m_1\cap Q_2 = \emptyset$ and $m_2\cap Q_1 =\emptyset$ (here, $\cap$ denotes the set theoretic intersection).
    Since $P$ is strongly connected, there exists a sequence $\sigma_0,\sigma_1, \ldots, \sigma_k$ of maximal chains such that $m_1=\sigma_0$, $m_2=\sigma_k$, and $\sigma_j\cap\sigma_{j+1}$ is a chain of length $\rank(P)-1$ for all $0\leq j\leq k-1$.
    Let $l$ be the smallest integer such that $\sigma_l\cap Q_2 \neq \emptyset$. 
    Since $\sigma_{l-1}\cap Q_2 = \emptyset$ by the choice of $l$ and $\sigma_{l-1}\cap \sigma_l$ is a chain of length $\rank(P)-1$, we get that $\sigma_l \cap Q_2$ is a singleton, say $\{a\}$. 
    First, assume that the height of $a$ is $i+1$ in $P$.
    Let $b\in \sigma_l$ be such that $b \lessdot a$. 
    Since $P$ is pure, height of $b$ is $i$ in $P$.
    Also, $b\notin Q_2$ because $\sigma_l \cap Q_2=\{a\}$; thus $b\in Q_1$
    which is a contradiction.
    Similar argument follows when height of $a$ is $i$ in $P$. 
    This completes the proof.
\end{proof}

Now the aim of this section is to prove \eqref{thm:cmpermutation:poset} $\implies$ \eqref{thm:cmpermutation:shellable} of the Theorem~\ref{thm:cmpermutation}. 
Let $\tau$ be a permutation on $[n]$. 
$\tau$ gives a linear order on $[n]$ as follows: for $i,j\in [n]$, $i<j$ in $\tau$ if there exist $a, b\in [n]$ with $a<b$ in $\naturals$ such that $\tau(a) =i$ and $\tau(b) =j$.
By abuse of terminology, we say that the permutation $\tau$ is a linear order.
We write $\tau$ as $[\tau_1,\ldots, \tau_n]$ where $\tau_a \coloneqq \tau(a)$ for all $a\in [n]$.
Note that $\tau_a <\tau_b$ in $\tau$ if and only if $\tau_b$ is on the right side of $\tau_a$ in $\tau$ for any $a,b\in [n]$.
For two permutations $\sigma$ and $\tau$, $P_{\sigma,\tau}$ denotes the poset that is the intersection of $\sigma$ and $\tau$. 
We start with an observation that a dimension two poset is isomorphic to a poset that is an intersection of the identity permutation and an another permutation.

\begin{proposition}
    \label{prop:perm_graph}
    Let $\sigma$ and $\tau$ be two permutations on $[n]$. 
Then, there exists a permutation $\pi$ such that $P_{\sigma,\tau}\iso P_{\id,\pi}$, where $\id$ is the identity permutation. 
\end{proposition}

\begin{proof}
    
Let $\pi = {\sigma}^{-1}\tau$. Define a map $\varphi : P_{\id,\pi} \to P_{\sigma,\tau}$ by $\varphi(j) = \sigma(j)$. 
Clearly, $\varphi$ is well-defined and it is a bijection.
It suffices to show that $i< j$ in $P_{\id,\pi}$ if and only if $\sigma(i)< \sigma(j)$ in $P_{\sigma,\tau}$.

Suppose that $i< j$ in $P_{\id,\pi}$, i.e., $i<j$ in $\naturals$ and there exist $a, b\in [n]$ with $a<b$ such that $\pi(a) =i$ and $\pi(b) =j$.
Clearly, $\sigma(i) <\sigma(j)$ in $\sigma$.
Also, we have $\pi(a) = {\sigma}^{-1}\tau(a) = i$; thus $\tau(a)=\sigma(i)$. Similarly, $\tau(b)=\sigma(j)$. 
So, $\tau(a)= \sigma(i)<\sigma(j)=\tau(b)$ in $\tau$. 
Hence $\sigma(i) <\sigma(j)$ in $P_{\sigma,\tau}$.

Conversely, suppose that $\sigma(i)< \sigma(j)$ in $P_{\sigma,\tau}$, i.e., $i<j$ in $\naturals$ and there exist $a, b\in [n]$ with $a<b$ such that $\tau(a) =\sigma(i)$ and $\tau(b) =\sigma(j)$.
Therefore, ${\sigma}^{-1}\tau(a) =i$ and ${\sigma}^{-1}\tau(b)=j$. So $i<j$ in $\pi$. Hence, $i< j$ in $P_{\id,\pi}$.
\end{proof}

\begin{examplebox}
    We illustrate the above proposition now. 
    Let $\sigma=[2, 3, 1, 4, 5]$ and $\tau =[3, 2, 1, 5, 4]$ be two permutations.
    Then, $P_{\sigma,\tau}$ is as shown in Figure~\ref{fig:poset}.
    Let $\pi = {\sigma}^{-1}\tau$. 
    Note that $\pi = [2, 1, 3, 5, 4]$ and $P_{\id,\pi}$ is as shown in Figure~\ref{fig:posetidentity}. 
    Also, it is immediate that $j \mapsto \sigma(j)$ is an isomorphism from $P_{\id,\pi}$ to $P_{\sigma,\tau}$.
\end{examplebox}

\begin{figure}[h]
	\begin{subfigure}[t]{7cm}
        \centering
		\begin{tikzpicture}[scale=1]
            \draw[] (0,0)--(.5,1)--(0,2)
			(1,0)--(.5,1)--(1,2);
            \filldraw[black] (1,2) circle (1pt) node[anchor=south]  {$5$};
            \filldraw[black] (1,0) circle (1pt) node[anchor=north]  {$3$};
            \filldraw[black] (0,2) circle (1pt) node[anchor=south]  {$4$};
            \filldraw[black] (.5,1) circle (1pt) node[anchor=west]  {$1$};
            \filldraw[black] (0,0) circle (1pt) node[anchor=north]  {$2$};
		\end{tikzpicture}
        \caption{$P_{\sigma,\tau}$}\label{fig:poset}
        \end{subfigure}
	\quad
	\begin{subfigure}[t]{7cm}
        \centering
		\begin{tikzpicture}[scale=1]
            \draw[] (0,0)--(.5,1)--(0,2)
			(1,0)--(.5,1)--(1,2);
            \filldraw[black] (1,2) circle (1pt) node[anchor=south]  {$5$};
            \filldraw[black] (1,0) circle (1pt) node[anchor=north]  {$2$};
            \filldraw[black] (0,2) circle (1pt) node[anchor=south]  {$4$};
            \filldraw[black] (.5,1) circle (1pt) node[anchor=west]  {$3$};
            \filldraw[black] (0,0) circle (1pt) node[anchor=north]  {$1$};
	\end{tikzpicture}
    \caption{$P_{\id,\pi}$}\label{fig:posetidentity}
        \end{subfigure}
        \caption{}   
\end{figure}

\begin{definitionbox}\label{def:posetdim2}
    Let $P$ be a poset of dimension two. 
    By relabeling, we may assume that $P$ is a poset on the set $[n]$. 
    So $P$ can be written as an intersection of two permutations, say $\sigma$ and $\tau$. 
    By Proposition~\ref{prop:perm_graph}, there exists a permutation $\pi$ such that $P \iso P_{\id,\pi}$, where $P_{\id,\pi}$ is the poset that is the intersection of the identity permutation and $\pi$.
    We denote $P_{\id,\pi}$ by $P_\pi$.
\end{definitionbox}

The idea of the proof of \eqref{thm:cmpermutation:poset} $\implies$ \eqref{thm:cmpermutation:shellable} of Theorem~\ref{thm:cmpermutation} is motivated by~\cite{Du03ShellablePlanarPosets}. 
Let $P_\pi$ be a poset as defined in Definition~\ref{def:posetdim2}.
For $0\leq i\leq \rank(P_\pi)$, let $P_i$ be the set of all height $i$ elements of $P_\pi$.
For all $i$, define a linear order $<_i$ on $P_i$ as following:
\[x<_i y \ \text{if and only if} \ x>y \ \text{in} \ \naturals.\]

For $x\in P_\pi$, let $U(x)$ be the set of all elements of $P_\pi$ that covers $x$. 
If $x\in P_i$, then $y\in P_{i+1}$ for all $y\in U(x)$ if $P_\pi$ is pure. 
For $x\in P_\pi$, define $x_{\min}:= \min(U(x))$ and $x_{\max}:=\max(U(x))$.

We make few observations which directly follows from the definition of $P_\pi$ and of the linear order $<_i$.

\begin{observationbox}\label{obs:permutation}
    Let $P_\pi$ be a pure poset. We have
    \begin{asparaenum}
    \item\label{obs:antichain}
        If $x<_i y$ in $P_i$, then $\pi$ has the form $[\ldots, x, \ldots, y,\ldots]$ because $x$ and $y$ are incomparable in $P_\pi$ and $x> y$ in $\naturals$.
     
    \item\label{obs:cover}
        If $y \in U(x)$ for some $x,y\in P_\pi$, then $x<y$ in $\naturals$. 
        Also, $y$ is on the right side of $x$ in $\pi$.
    
    \item\label{obs:threeantichain}
          If $x_{\min} <_i y <_i x_{\max}$ for some $x\in P_{i-1}$ and $y\in P_i$, then $y\in U(x)$.
          In fact, by~\eqref{obs:antichain} and~\eqref{obs:cover}$, \pi$ has the form $[\ldots, x, \ldots, x_{\min}, \ldots, y,\ldots, x_{\max},\ldots]$. 
    Also, note that $x<x_{\max}$ in $\naturals$ because $x_{\max}\in U(x)$; thus $x<y$ in $\naturals$.  
    Hence, $y\in U(x)$.
    \end{asparaenum}
\end{observationbox}

    \begin{Lemma}\label{lemma:maximalchaininterval}
    Let $P_\pi$ be a poset that satisfies~\eqref{thm:cmpermutation:poset} of the Theorem~\ref{thm:cmpermutation}. 
    Let $[x,y]$ be an interval in $P_\pi$ such that $\height(x)=i,$ $\height(y)=j$ and $j-i\geq 2$. 
    Let $x=y_i\lessdot y_{i+1}\lessdot \cdots\lessdot y_j=y$ be a chain in $[x,y]$ such that for all $k$ with $i<k<j$ there exists an $x_k\in [x,y]$, $x_k\lessdot_k y_k$.
    Then, there exists some integer $k'$, $i<k'<j$ such that $y_{k'-1}\lessdot x_{k'}\lessdot y_{k'+1}$.
\end{Lemma}

\begin{proof}
    First, we claim the following: if $x \lessdot_l y$ in $P_l$ for some $x,y\in P_l$ and $0\leq l\leq \rank(P_\pi)-1$, then $y_{\min} \leq_{l+1} x_{\max}$ in $P_{l+1}$.

    Assume that the claim holds.
    If $j-i=2$, then we can take $k'=i+1$.
    Now assume that $j-i>2$. Consider $x_{i+1}$. If $x_{i+1}\lessdot y_{i+2}$, then we can take $k'=i+1$.
    Otherwise, if $x_{i+1}\nless y_{i+2}$, then $y_{i+1}\lessdot x_{i+2}$ by the claim.
    Now, consider $x_{i+2}$. 
    If $x_{i+2} \lessdot y_{i+3}$, then we can take $k'=i+2$. 
    Otherwise if $x_{i+2}\nless y_{i+3}$, then $y_{i+2}\lessdot x_{i+3}$ by the claim.
    Proceeding in this way and using the case $j-i=2$, we find the desired $x_{k'}$,
    which completes the proof.

    We now prove the claim.
    On the contrary, suppose that there exists a $l\in \{0,1,\ldots,\rank(P_\pi)-1\}$ and $x,y\in P_l$ such that $x \lessdot_l y$ in $P_l$ and $x_{\max} <_{l+1} y_{\min}$ in $P_{l+1}$. 
    We show that the induced subposet $Q$ of $P_\pi$ consisting $P_l$ and $P_{l+1}$ is disconnected.
    
    Define \[Q_1 = \{p\in P_\pi :\text{either}\ p\in P_l \ \text{and}\ p\leq_l x \ \text{or} \  p\in P_{l+1} \ \text{and}\ p<_{l+1} y_{\min}\}\] 
    and 
    \[Q_2 = \{p\in P_\pi :\text{either}\ p\in P_l \ \text{and}\ y\leq_l p \ \text{or} \  p\in P_{l+1} \ \text{and}\ y_{\min}\leq_{l+1} p\}.\]
    
    We show that $Q$ is the disjoint union of the subposets $Q_1$ and $Q_2$. 
    Suppose that there exists an edge between $Q_1$ and $Q_2$. 
    So, either there exists an $x' <_{l} x$ with $y_{\min}\leq_{l+1}x'_{\max}$ or there exists an $y'\in P_l$ with $y<_l y'$ and $y'_{\min}<_{l+1}y_{\min}$.
    We consider both cases separately:

    $(i)$ There exists an $x' <_{l} x$ with $y_{\min}\leq_{l+1}x'_{\max}$. 
    Then, $x<x',$ $x'_{\max}\leq y_{\min}$ in $\naturals$.
    Using Observation~\ref{obs:permutation}, we get that $\pi$ has the form $[\ldots, x', \ldots, x, \ldots, x_{\max},\ldots,y_{\min},\ldots,x'_{\max},\ldots]$,
    in fact $x$ is on the right side of $x'$ in $\pi$ by~\eqref{obs:antichain} of the Observation~\ref{obs:permutation}, $x_{\max}$ is on the right side of $x$ in $\pi$ by~\eqref{obs:cover} of the Observation~\ref{obs:permutation}, and the position of $x_{\max},y_{\min},x'_{\max}$ is by~\eqref{obs:antichain} of the Observation~\ref{obs:permutation}.
    By \eqref{obs:cover} of the Observation~\ref{obs:permutation}, $x' <x'_{\max}$ in $\naturals$. 
    Since $x<x'<x'_{\max}$ in $\naturals$ and $x'_{\max}$ is on the right side of $x$ in $\pi$, we get that $x'_{\max}\in U(x)$. Which is a contradiction.

    $(ii)$  There exists an $y'\in P_l$ with $y<_l y'$ and $y'_{\min}<_{l+1}y_{\min}$. 
    Then, $y'<y,$ $y_{\min}<y'_{\min}$ in $\naturals$.
    Using Observation~\ref{obs:permutation}, $\pi$ has the form $[\ldots, y, \ldots, y', \ldots, y'_{\min},\ldots,y_{\min},\ldots]$. 
    Note that $y'_{\min}$ is on the right side of $y$ in $\pi$. 
    Also, note that $y<y'_{\min}$ in $\naturals$ because $y\leq y_{\min} \leq y'_{\min}$ in $\naturals$.
    Hence $y'_{\min}\in U(y)$ which is a contradiction. This completes the proof of the claim. 
\end{proof}

We are now ready to prove \eqref{thm:cmpermutation:poset} $\implies$ \eqref{thm:cmpermutation:shellable} of our main theorem. 

\begin{proof}[Proof of \eqref{thm:cmpermutation:poset} $\implies$ \eqref{thm:cmpermutation:shellable}]
    By Definition~\ref{def:posetdim2}, it suffices to prove the result for the posets $P_\pi$, where $P_\pi$ is as defined in Definition~\ref{def:posetdim2}.
    Assume that $P_\pi$ is a poset of rank $r$ that satisfies~\eqref{thm:cmpermutation:poset}.
    If $\pi = [n,n-1,\ldots,1]$, then $P_\pi$ is an antichain. 
    So antichains have dimension two.
    It follows from the definition of the shellability that antichains are shellable.  
    So we may assume that $P_\pi$ is not an antichain. 
    
    Consider the permutation $\pi' = [0,\pi,n+1]$ on $n+2$ elements. 
    Then, $P_{\pi'} = P_\pi\cup \{0,n+1\}$, where $0$ and $n+1$ are the minimal and the maximal elements of $P_{\pi'}$ respectively. 
    Note that $P_{\pi'}$ is a pure poset of rank $r+2$ and it satisfies the hypothesis of~\eqref{thm:cmpermutation:poset}.
    Since every interval of a shellable poset is also shellable~\cite[Proposition\ 8.2]{BW83lexicographicallyShellablePosets}, we may replace $\pi$ by $\pi'$.

    Let $C : x_0\lessdot x_1\lessdot \cdots \lessdot x_{r+2}$ and $C' : y_0\lessdot y_1\lessdot \cdots \lessdot y_{r+2}$ be two maximal chains of $P_\pi$.
    Note that $x_0=y_0 = 0$ and $x_{r+2}=y_{r+2} =n+1$.
    Let $j' = \max\{i\in [r+2] : x_i \neq y_i\}$. Define a linear order $<_E$ on the maximal chains of $P_\pi$ as follows:
    \[C<_E C' \ \text{if and only if} \ x_{j'} <_{j'} y_{j'}.\]
    Under the above notations,  assume that $C<_E C'$.
    Let $i' = \max\{i<j' : x_i =y_i\}.$ 
    Then, $x_{i'} =y_{i'}\lessdot y_{i'+1}\lessdot \cdots \lessdot y_{j'+1}=x_{j'+1}$ is a maximal chain in $[x_{i'},y_{j'+1}]$. 
    We proceed in the following cases:

    \begin{asparaenum}
        \item\label{case} 
            If $x_k<_k y_k$ for all $i'<k<j'+1$, then the maximal chain $x_{i'} =y_{i'}\lessdot y_{i'+1}\lessdot \cdots \lessdot y_{j'+1}=x_{j'+1}$ in the interval $[x_{i'},y_{j'+1}]$ satisfies the hypothesis of Lemma~\ref{lemma:maximalchaininterval}, 
            i.e., for all $k$, $i'<k<j'+1$ there exists a $z_k\in [x_{i'},y_{j'+1}]$ such that $z_k\lessdot_k y_k$ because $<_k$ is linear order and $x_k<_k y_k$ for all $k$.
            So there exists a $k'$, $i'<k'<j'+1$ such that $y_{k'-1}\lessdot z_{k'}\lessdot y_{k'+1}$ by~Lemma~\ref{lemma:maximalchaininterval}. 
            If we let \[C'' : y_0\lessdot y_1\lessdot \cdots \lessdot y_{k'-1}\lessdot z_{k'} \lessdot y_{k'+1} \lessdot \cdots \lessdot y_{r+2},\]
            then $C''<_E C',$ $y_{k'} \in C'\setminus C$ and   $C'\setminus C'' = \{y_{k'}\}$.
    
        \item
            There exists a $k$ with $i'<k<j'+1$ such that $y_k<_k x_k$. 
            Let $l = \max\{k : i'<k<j'+1 \ \text{and}\ y_k<_k x_k\}$.
            First, we show that $y_l < x_{l+1}$ in $P_\pi$. 
            By the choice of $l$, we have $x_{l+1}<_{l+1} y_{l+1}$;
            so $y_{l+1}<x_{l+1}$ in $\naturals$. Also, $y_l<y_{l+1}$ in $P_\pi$. Thus $y_{l}< y_{l+1}<x_{l+1}$ in $\naturals$.
            Under the given conditions,  $\pi$ has the form $[\ldots, y_l, \ldots, x_l, \ldots, x_{l+1},\ldots,y_{l+1},\ldots]$.
            Observe that $x_{l+1}$ is on the right side of $y_l$ in $\pi$.
            Therefore, $y_l<x_{l+1}$ in $P_\pi$.
            
            Now consider the interval $[y_l,y_{j'+1}]$. Note that $j'+1-l\geq 2$ because $x_{j'} <_{j'} y_{j'}$.
            Also, observe that for all $l<k<j'+1$, $x_k \in [y_l,y_{j'+1}]$ because $y_l<x_{l+1}$ in $P_\pi$, and $x_k<_k y_k$ by the choice of $l$.
            Thus for all $k$, $l<k<j'+1$ there exists a $z_k\in [y_{l},y_{j'+1}]$ such that $z_k\lessdot_k y_k$.
            Thus, the maximal chain $y_{l}\lessdot y_{l+1}\lessdot \cdots \lessdot y_{j'+1}=x_{j'+1}$ in the interval $[y_l,y_{j'+1}]$ satisfies the hypothesis of Lemma~\ref{lemma:maximalchaininterval}.
            So there exists a $k'$, $l<k'<j'+1$ such that $y_{k'-1}\lessdot z_{k'}\lessdot y_{k'+1}$. 
            Therefore, we can repeat the argument of \eqref{case} to complete the proof. 
    \end{asparaenum}
    Therefore, $<_E$ is a shelling order on the maximal chains of $P_\pi$. This completes the proof.
\end{proof}

We see that the Theorem~\ref{thm:cmpermutation} helps us to characterize the Cohen-Macaulay permutation graphs.
Let $G$ be a graph on $[n]$.
Let $S = K[x_1,\ldots,x_n]$ be a polynomial ring over a field $K$.
Let $I_G = (x_ix_j : \{i,j\} \ \text{is an edge of} \ G)$ be the {\em edge ideal} of $G$.
We say that $G$ is {\em Cohen-Macaulay} if $S/I_G$ is Cohen-Macaulay.

Now assume that $G$ is a permutation graph.
Then $G$ is the intersection graph of the permutation diagram consisting of horizontal lines $l_0$ and $l_1$.
Assume that $l_0$ and $l_1$ are labeled by the permutations $\pi_0$ and $\pi_1$ respectively.
Let $P$ be the poset that is the intersection of $\pi_0$ and $\pi_1$. Then $G$ is the co-comparability graph of $P$ by~\cite[Theorem\ 1]{GRU83comparabilitygraphs}, and the dimension of $P$ is at most two.
The Stanley-Reisner ideal of the order complex of $P$ coincide with the edge ideal of $G$.
Thus, by~\cite[Theorem\ 1]{ReisnerCMrings}, we get that $P$ is Cohen-Macaulay if and only if $G$ is Cohen-Macaulay. 
If the dimension of $P$ is one, then $P$ is a linear order; thus, $I_G$ is the trivial ideal. 
Hence $G$ is Cohen-Macaulay. 
When the dimension of $P$ is two, we can use Theorem~\ref{thm:cmpermutation} to check whether $G$ is Cohen-Macaulay. 
Consequently, $G$ is Cohen-Macaulay over any field.

\begin{figure}[h]
\begin{center}
    \begin{tikzpicture}[scale=1.3]
\draw[] (0,0)--(0,1)--(1,0)--(1,1)--(0,0) 
(2,0)--(2,1)--(0,0)
(0,1)--(.5,2)--(2,1)
(1,1)--(1.5,2)--(2,1)
;

\filldraw[black] (0,0) circle (.5pt) node[anchor=north]  {3};
\filldraw[black] (1,0) circle (.5pt) node[anchor=north] {2};
\filldraw[black] (2,0) circle (.5pt) node[anchor=north] {1};
\filldraw[black] (0,1) circle (.5pt) node[anchor=east] {4};
\filldraw[black] (1,1) circle (.5pt) node[anchor=west]  {5};
\filldraw[black] (2,1) circle (.5pt) node[anchor=west] {6};
\filldraw[black] (.5,2) circle (.5pt) node[anchor=south] {7};
\filldraw[black] (1.5,2) circle (.5pt) node[anchor=south] {8};
\end{tikzpicture}
\caption{}\label{figure:non-cm}
\end{center}
\end{figure}

It follows from the \cite[Proposition\ 11.7]{Bjornertopologicalmethods} and Lemma~\ref{lem:stronglyconnected} that the \eqref{thm:cmpermutation:poset} of Theorem~\ref{thm:cmpermutation} is a necessary condition for a poset to be Cohen-Macaulay.
    In the following example, we show that it is not a sufficient condition.
    More precisely, we show that \eqref{thm:cmpermutation:poset} $\implies$ \eqref{thm:cmpermutation:cm} and  \eqref{thm:cmpermutation:poset} $\implies$ \eqref{thm:cmpermutation:shellable} of the Theorem~\ref{thm:cmpermutation} may not be true when the dimension of the poset is $\geq 3$.

\begin{examplebox}
    Consider the poset $P$ as shown in Figure~\ref{figure:non-cm}.
    A SageMath~\cite{[sagemath]} computation shows that the dimension of $P$ is three. 
    In fact, $P$ is the intersection of the following permutations: 
    $[1, 3, 6, 2, 4, 7, 5, 8],$ $[2, 3, 4, 5, 1, 6, 7, 8]$ and $[3, 1, 2, 6, 5, 8, 4, 7]$.
    Clearly, $P$ satisfies the hypothesis of \eqref{thm:cmpermutation:poset} of the Theorem~\ref{thm:cmpermutation}.
    Note that $\link(\Delta(P),\{2\})=\{\{4,7\},\{5,8\}\}$ which is disconnected. 
    Since $\dim(\link(\Delta(P),\{2\}))=1$ and $\widetilde{H}_0(\link(\Delta(P),\{2\}), K)= K$, $P$ is not Cohen-Macaulay. 
    Hence, $P$ is not shellable. 
\end{examplebox}
    
All examples we have computed suggest that shellability and Cohen-Macaulayness coincide for dimension three posets.
    This, along with the Theorem~\ref{thm:cmpermutation}, motivates us to ask the following question: 
    What is the least $d\in \naturals$ such that there exists a Cohen-Macaulay poset of dimension $d$ that is not shellable?

%\bibliography{../ref.bib}{}
%\bibliographystyle{alpha}   

\end{document}